\documentclass[11pt]{amsart}

\usepackage{geometry}
\usepackage{todonotes}
\usepackage{amsfonts, adjustbox, amssymb, amscd}
\numberwithin{equation}{section}

\usepackage{bm}
\usepackage{verbatim}
\usepackage{mathrsfs}
\usepackage{graphicx}
\usepackage{tikz-cd}
\usepackage{subcaption}
\usepackage{listings}
\usepackage{subfiles}
\usepackage[toc,page]{appendix}
\usepackage{mathtools}
\usepackage{comment}
\usepackage{enumerate}
\usepackage{enumitem}
\usepackage[all]{xy}

\usepackage{graphicx}
\usepackage{appendix}
\usepackage{hyperref}
\hypersetup{
    colorlinks=true,
    citecolor=red,
    linkcolor=blue,
    filecolor=magenta,      
    urlcolor=red,
}
\newcommand{\bb}{\bm{b}}
\newcommand{\Mm}{{\bf{M}}}

\newcommand{\Qq}{\mathbb{Q}}

\newcommand{\Rr}{\mathbb{R}}

\newcommand{\Supp}{\operatorname{Supp}}

\newcommand{\Nlc}{\operatorname{Nlc}}


\newcounter{parentnumber}

\newtheorem{thm}{Theorem}[section]

\newtheorem{lem}[thm]{Lemma}

\theoremstyle{definition}
\newtheorem{defn}[thm]{Definition}
\newtheorem{ques}[thm]{Question}
\theoremstyle{definition}
\newtheorem{rem}[thm]{Remark}

\theoremstyle{definition}

\begin{document}

%\title{On equivalence of Prokhorov-Shokurov and Li's effective adjunction conjectures}
\title[A generalized non-vanishing theorem on surfaces]{A generalized non-vanishing theorem on surfaces}
\author{Jihao Liu, Lingyao Xie}

\date{\today}

\begin{abstract}
We show that the anti-canonical bundle of any $\Qq$-factorial surface is numerically effective if and only if it is pseudo-effective. To prove this, we establish a numerical non-vanishing theorem for surfaces polarized with pseudo-effective divisors. The latter answers a question of C. Fontanari.
\end{abstract}

\address{Department of Mathematics, Peking University, No. 5 Yiheyuan Road, Haidian District, Peking 100871, China}
\email{liujihao@math.pku.edu.cn}

\address{Department of Mathematics, University of California, San Diego, 9500 Gilman Drive \# 0112, La Jolla, CA 92093-0112, USA}
\email{l6xie@ucsd.edu}

\maketitle

\tableofcontents

\section{Introduction}\label{sec:Introduction}

We work over the field of complex numbers $\mathbb C$. Without further remarks, all pairs and generalized pairs are assumed to have $\Qq$-coefficients in this paper.

The \emph{non-vanishing} conjecture predicts that any the log canonical $\Qq$-divisor $K_X+B$ of any projective lc pair $(X,B)$ is effective if it is pseudo-effective (cf. \cite[Conjecture 2.1]{BCHM10}). This is a very difficult conjecture and is only known up to dimension three \cite[1.2 Theorem]{KMM94}. A variation of the non-vanishing conjecture, the so-called \emph{generalized non-vanishing} conjecture, predicts that the generalized log canonical $\Qq$-divisor $K_X+B+\Mm_X$ of any projective lc generalized pair $(X,B,\Mm)$ is numerically effective (i.e. $K_X+B+\Mm_X\equiv D\geq 0$ for some $\Qq$-divisor $D$) if it is pseudo-effective  \cite[Conjecture 1.2]{HL20} (see \cite[Question 3.5]{BH14}, \cite{LP20a} for the case when $\Mm$ descends to $X$). This conjecture is important for the study of generalized pairs, which was introduced by Birkar and Zhang \cite{BZ16} and plays a crucial role in modern birational geometry. For other works dedicated to the generalized non-vanishing conjecture, we refer the reader to \cite{LP20b,BS23,LMPTX23,Mul23}.

C. Fontanari recently asked the following question on whether the generalized non-vanishing holds in greater generality. That is, instead of the polarization of a nef ($\bb$-)$\Qq$-divisor, we only polarize with a pseudo-effective $\Qq$-divisor:
\begin{ques}[{\cite[Question 2]{Fon24}}]\label{ques: Fon24 q2}
    Let $(X,B)$ be a $\Qq$-factorial projective lc $\Qq$-pair and $M$ a pseudo-effective $\mathbb Q$-divisor on $X$ such that $K_X+B+M$ is pseudo-effective. Is there an effective $\Qq$-divisor $D$ on $X$ such that $K_X+B+M\equiv D$?
\end{ques}

\cite{Fon24} pointed out that Question \ref{ques: Fon24 q2} remained open in general even for surfaces, where some special cases were proven in \cite[Propositions 3,5, Theorems 4,6]{Fon24}. 

In this paper, we answer Fontanari's question for surfaces in full generality. In fact, we do not even need $X$ to be lc:

\begin{thm}\label{thm: main theorem}
    Let $X$ be a $\Qq$-factorial projective surface and $L$ a pseudo-effective $\Qq$-divisor on $X$ such that $K_X+L$ is pseudo-effective. Then $K_X+L$ is numerically effective.
\end{thm}

As a consequence, we immediately deduce the following result:

\begin{thm}\label{thm: anti pe is ne}
      Let $X$ be a $\Qq$-factorial projective surface such that $-K_X$ is pseudo-effective. Then $-K_X$ is numerically effective.  
\end{thm}

\noindent\textbf{Acknowledgement}. The authors would like to thank Jingjun Han for useful discussions. They would like to thank the referee for many useful comments. The second author had been partially supported by NSF research grants no. DMS-1801851 and DMS-1952522, as well as a grant from the Simons Foundation (Award Number: 256202).

\section{Preliminaries}

We follow the standard notations and definitions as in \cite{BCHM10,KM98}. 

\begin{defn}
    Let $X$ be a normal projective variety and $D$ a $\Qq$-divisor on $X$. We say that $D$ is \emph{numerically effective} if $D-D'\equiv 0$ for some $\Qq$-divisor $D'\geq 0$ on $X$.
\end{defn}

\begin{lem}\label{lem: num eff r to q}
   Let $X$ be a normal projective variety and $D$ a $\Qq$-divisor on $X$. Then $D$ is numerically effective if and only if  $D-D'\equiv 0$ for some $\Rr$-divisor $D'\geq 0$ on $X$.
\end{lem}
\begin{proof}
    We may write $D'=D_0+\sum_{i=1}^c r_iD_i$ where $D_0,\dots,D_c$ are $\Qq$-divisors and $1,r_1,\dots,r_c$ are linearly independent over $\mathbb Q$. Since $D'\geq 0$, there exists an open set $U\ni (r_1,\dots,r_c)$ in $\mathbb R^c$, such that $D_0+\sum_{i=1}^cv_iD_i\geq 0$ for any $(v_1,\dots,v_c)\in U$. Then 
    $$(D-D_0)+\sum_{i=1}^cr_iD_i\equiv 0.$$
    By \cite[Lemma 5.3]{HLS19}, $D_i\equiv 0$ for any $1\leq i\leq c$. We let $r_1'',\dots,r_c''$ be rational numbers such that $(r_1'',\dots,r_c'')\in U$ and let $D'':=D_0+\sum_{i=1}^cr_i''D_i$. Then $D\equiv D''\geq 0$ and $D''$ is a $\Qq$-divisor. The lemma follows.
\end{proof}

\begin{defn}
Let $(X,B)$ be a pair. A \emph{non-lc place} of $(X,B)$ is a prime divisor $D$ over $X$ such that the log discrepancy of $D$ with respect to $(X,B)$ is negative. A \emph{non-lc center} of $(X,B)$ is the center of a non-lc place of $(X,B)$ on $X$. The \emph{non-lc locus} $\Nlc(X,B)$ of $(X,B)$ is the union of all non-lc centers of $(X,B)$. 
\end{defn}

\begin{defn}
Let $X$ be a normal projective variety. We let $NE(X)$ be the cone of effective $1$-cycles on $X$ and let $\overline{NE}(X)$ be the closure of $NE(X)$. For any $\Rr$-Cartier $\Rr$-divisor $D$ on $X$, we define $\overline{NE}(X)_{D\geq 0}:=\{R\mid R\in\overline{NE}(X),D\cdot R\geq 0\}$. For any pair $(X,B)$, we define $\overline{NE}(X)_{\Nlc(X,B)}$ be the image of $\overline{NE}(\Nlc(X,B))\rightarrow\overline{NE}(X)$ where $\Nlc(X,B)$ is the non-lc locus of $(X,B)$. For any extremal face $F$ of $\overline{NE}(X)$, we say that $F$ is \emph{relatively ample at infinity with respect to} $(X,B)$  if $$F\cap\overline{NE}(X)_{\Nlc(X,B)}=\{0\}.$$
\end{defn}

\begin{thm}[Cone theorem for not necessarily lc pairs, {cf. \cite[Theorem 5.10]{Amb03}, \cite[Theorems 4.5.2, 6.7.4]{Fuj17}}]\label{thm: cone theorem for not necessarily lc pairs}
Let $(X,B)$ be a projective pair and let $\{R_j\}_{j\in\Lambda}$ be the set of $(K_X+B)$-negative extremal rays in $\overline{NE}(X)$ that are rational and relatively ample at infinity with respect to $(X,B)$. Then:
\begin{enumerate}
    \item $$\overline{NE}(X)=\overline{NE}(X)_{K_X+B\geq 0}+\overline{NE}(X)_{\Nlc(X,B)}+\sum_{j\in\Lambda} R_j.$$
    \item Each $R_j$ is spanned by a rational curve $C_j$ such that 
    $$0<-(K_X+B)\cdot C_j\leq 2\dim X.$$
\end{enumerate}
\end{thm}

\begin{thm}\label{thm: hl20 1.5}
    Let $X$ be a smooth projective surface and $P$ a nef $\Qq$-divisor on $X$ such that $K_X+P$ is pseudo-effective. Then $K_X+P$ is numerically effective.
\end{thm}
\begin{proof}
It follows from \cite[Theorem 1.5]{HL20} and Lemma \ref{lem: num eff r to q}.
\end{proof}

\begin{lem}\label{lem: num 0 pushforward 0}
    Let $h: X'\rightarrow X$ be a birational morphism such that $X$ is $\Qq$-factorial normal projective. Let $D$ be a $\Qq$-divisor on $X'$ such that $D\equiv 0$. Then $h_*D\equiv 0$.
\end{lem}
\begin{proof}
  Let $E:=D-h^*h_*D$. Then $E\equiv_X0$, $h_*E=0$, and $E$ is exceptional$/X$. By applying the negativity lemma twice, we have $E=0$. Thus $h^*h_*D\equiv 0$. By the projection formula, $h_*D\equiv 0$. 
\end{proof}

\begin{lem}\label{lem: smooth surface ne ray is -1 curve}
    Let $X$ be a smooth projective surface and $C$ a rational curve which spans a $K_X$-negative extremal ray, and the contraction $\pi: X\rightarrow X'$ of $C$ is a divisorial contraction. Then $C^2=-1$.
\end{lem}
\begin{proof}
Since $C$ is a rational curve, $(K_X+C)\cdot C=-2$. By the negativity lemma, we have $\pi^*K_{X'}+\lambda C=K_X$ for some positive real number $\lambda $, so $\lambda C\cdot C=K_X\cdot C<0$. Thus $K_X\cdot C<0$ and $C^2<0$. Since $X$ is smooth, $K_X\cdot C\leq -1$ and $C^2\leq -1$. Thus $C^2=-1$.
\end{proof}

\section{Induction approach to the main theorem}

The goal of this section is to prove the following theorem for any positive integer $\rho$:

\medskip

\noindent {\bf \hypertarget{thm: A}{Theorem A$_\rho$}} {\bf.} 
{\it    Let  $X$ be a smooth projective surface such that $\rho(X)\leq\rho$ and let $L$ be a pseudo-effective $\Qq$-divisor on $X$ such that $K_X+L$ is pseudo-effective. Then $K_X+L$ is numerically effective.}

\medskip

In this section, we shall prove Theorem \hyperlink{thm: A}{A} by induction on $\rho$.

\begin{lem}\label{lem: nef case}
 Let $\rho\geq 2$ be an integer. Assume that Theorem \hyperlink{thm: A}{A}$_{\rho-1}$ holds.

Let $X$ be a smooth projective surface such that $\rho(X)\leq\rho$, and let $L$ be a pseudo-effective $\Qq$-divisor on $X$ such that $K_X+L$ is nef. Then $K_X+L$ is numerically effective.
\end{lem}
\begin{proof}
    By induction hypothesis we may assume that $\rho(X)=\rho$. Let $P$ and $N$ be the positive and negative part of the Zariski decomposition of $L$ respectively. We may let $C_1,\dots,C_n$ be the irreducible components of $\Supp N$. Then $P\cdot C_i=0$ for each $i$. In the following, we apply induction on $n$. When $n=0$, the lemma follows from Theorem \ref{thm: hl20 1.5}. In the following, we may assume that $N\not=0$.

    Since the intersection matrix of $N$ is negative definite, there exists a unique $\Qq$-divisor $N'\geq 0$ on $X$ such that $\Supp N'=\Supp N$ and $N'\cdot C_i=-1$ for each $i$. Then there exists a unique positive rational number $c$ such that $N-cN'\geq 0$ and $\Supp(N-cN')\not=\Supp N$. We let $m$ be a positive integer such that $m(K_X+L)$ and $mcN'$ are Weil divisors. Since $X$ is smooth, $m(K_X+L)$ and $mcN'$ are Cartier. 
    
    Since $X$ is smooth,
    $$\Nlc(X,N-cN')\subset\Nlc(X,N)\subset\cup_{i=1}^nC_i.$$
    Therefore, by Theorem \ref{thm: cone theorem for not necessarily lc pairs},
    $$\overline{NE}(X)=\overline{NE}(X)_{K_X+N-cN'\geq 0}+\sum_{i=1}^nC_i+\sum_{j\in\Lambda}R_j,$$
    where each $R_j$ is a $(K_X+N-cN')$-negative extremal ray that is relative ample at infinity with respect to $(X,N-cN')$, and each $R_j$ is spanned by a rational curve $D_j$ such that
    $$0>(K_X+N-cN')\cdot D_j\geq -4.$$
Since $K_X+L$ is nef,
    $$(K_X+N-cN')\cdot C_i=(K_X+P+N-cN')\cdot C_i=(K_X+L)\cdot C_i-c(N'\cdot C_i)\geq c>0,$$
    for each $i$, so $[C_i]\in\overline{NE}(X)_{K_X+N-cN'\geq 0}$ for each $i$. Therefore,
    $$\overline{NE}(X)=\overline{NE}(X)_{K_X+N-cN'\geq 0}+\sum_{j\in\Lambda}R_j,$$
    and $D_j\not\subset\Supp N$ for any $j$. Since $P$ is nef, $\overline{NE}(X)_{K_X+N-cN'\geq 0}\subset\overline{NE}(X)_{K_X+P+N-cN'\geq 0}$. We let
    $$\Lambda_0:=\{j\in\Lambda\mid (K_X+P+N-cN')\cdot R_j<0\},$$
then
    $$\overline{NE}(X)=\overline{NE}(X)_{K_X+P+N-cN'\geq 0}+\sum_{j\in\Lambda_0}R_j.$$  
    
    If $\Lambda_0=\emptyset$, then $K_X+P+(N-cN')$ is nef. Let $L':=P+N-cN'$, then $P$ and $N-cN'$ are the positive and negative parts of the Zariski decomposition of $L'$ respectively, and $\Supp(N-cN')$ had $\leq n-1$ components. By induction on $n$, $K_X+L'$ is numerically effective. Note that here we use induction on $n$, the number of irreducible components of the negative part of the Zariski decomposition of $L$, rather than an induction on $\rho$. Thus $K_X+L=K_X+L'+cN'$ is numerically effective. Therefore, we may assume that $\Lambda_0\not=\emptyset$.

    For any $j\in\Lambda_0$, we have
    $$0>(K_X+P+N-cN')\cdot D_j\geq (K_X+N-cN')\cdot D_j\geq -4$$
    and
    $$(K_X+P+N)\cdot D_j\geq 0.$$
    Therefore, there exists a unique rational number $\lambda_j\in [0,1)$ such that
    $$(K_X+P+N-\lambda_jcN')\cdot D_j=0.$$
    We have
    $$\lambda_j=\frac{m(K_X+L)\cdot D_j}{m(K_X+L)\cdot D_j-m(K_X+P+N-cN')\cdot D_j}\in\left\{\frac{a}{a+b}\mid a,b\in\mathbb N, 0<b\leq 4m\right\}$$
    Since the set $\left\{\frac{a}{a+b}\mid a,b\in\mathbb N, 0<b\leq 4m\right\}\subset [0,1)$ satisfies the DCC condition (i.e. any descending sequence in the set stabilizes), there exists $j_0\in\Lambda_0$ such that $\lambda_{j_0}=\min\{\lambda_j\mid j\in\Lambda_0\}$. We let $L':=L-\lambda_{j_0}cN'=P+N-\lambda_{j_0}cN'$.
    
    By construction, $(K_X+L')\cdot D_j\geq 0$ for any $j\in\Lambda_0$, so $(K_X+L')\cdot R\geq 0$ for any $R\not\subset\overline{NE}(X)_{K_X+P+N-cN'\geq 0}$. Since $K_X+L$ is nef, for any $R\in\overline{NE}(X)_{K_X+P+N-cN'\geq 0}$, we have
    $$(K_X+L')\cdot R=(1-\lambda_{j_0})(K_X+L)\cdot R+\lambda_{j_0}(K_X+P+N-cN')\cdot R\geq 0.$$
  Since $\overline{NE}(X)=\overline{NE}(X)_{K_X+P+N-cN'\geq 0}+\sum_{j\in\Lambda_0}R_j$,  $K_X+L'$ is nef, and $(K_X+L')\cdot R_{j_0}=0$.

    Since $D_{j_0}\not\subset\Supp N$ and $(K_X+P+N-cN')\cdot D_{j_0}<0$, $K_X\cdot D_{j_0}<0$. Since $D_{j_0}$ spans the extremal ray $R_{j_0}$, $D_{j_0}$ is a $K_X$-negative extremal ray, and there exists a contraction $\pi: X\rightarrow \bar X$ of $D_{j_0}$. Since $\rho\geq 2$ and $X$ is smooth, by Castelnuovo's theorem, $\bar X$ is either a smooth surface or a curve. More precisely: If $\pi$ is not a divisorial contraction, then either $\bar X$ is a curve or $\bar X$ is a point, but since $\rho(X/\bar X)=1$ and $\rho\geq 2$, we have that $\bar X$ is a curve. If $\pi$ is a divisorial contraction, then since $D_{j_0}$ is a rational curve, by Lemma \ref{lem: smooth surface ne ray is -1 curve}, $D_{j_0}^2=-1$. By Castelnuovo's theorem, $\bar X$ is smooth.

    If $\bar X$ is a curve, then by the contraction theorem,
    $$K_X+L'\sim_{\mathbb Q,\pi}0.$$
    Since $K_X+L'$ is nef, there exists a nef $\Qq$-divisor $H$ on $\bar X$ such that $K_X+L'\sim_{\mathbb Q}\pi^*H$. Since $H$ is numerically effective, $K_X+L'$ is numerically effective, hence $K_X+L=K_X+L'+\lambda_{j_0}cN'$ is numerically effective, and we are done.

    If $\bar X$ is a smooth surface, then we let $\bar L':=\pi_*L$. By the contraction theorem (cf. \cite[Theorem 3.7(4)]{KM98}), 
    $$K_X+L'=\pi^*(K_{\bar X}+\bar L').$$
    Since $K_X+L'$ is nef, $K_{\bar X}+\bar L'$ is nef. By Theorem \hyperlink{thm: A}{A}$_{\rho-1}$, $K_{\bar X}+\bar L'$ is numerically effective, so $K_X+L'$ is numerically effective. Thus $K_X+L=K_X+L'+\lambda_{j_0}cN'$ is numerically effective, and we are done.
\end{proof}

\begin{lem}\label{lem: nef on n case}
 Let $\rho\geq 2$ be an integer. Assume that Theorem \hyperlink{thm: A}{A}$_{\rho-1}$ holds.

Let $X$ a smooth projective surface such that $\rho(X)\leq\rho$. Let $L$ be a pseudo-effective $\Qq$-divisor on $X$ such that $K_X+L$ is pseudo-effective. Let $P$ and $N$ be the positive part and negative part of the Zariski decomposition of $L$ respectively. Assume that $(K_X+L)\cdot C\geq 0$ for any irreducible component $C$ of $\Supp N$.

Then $K_X+L$ is numerically effective.
\end{lem}
\begin{proof}
 By induction hypothesis we may assume that $\rho(X)=\rho$. Let $C_1,\dots,C_n$ be the irreducible components of $\Supp N$. 
  Since $X$ is smooth,
    $$\Nlc(X,N-cN')\subset\Nlc(X,N)\subset\cup_{i=1}^nC_i.$$
    Therefore, by Theorem \ref{thm: cone theorem for not necessarily lc pairs},
    $$\overline{NE}(X)=\overline{NE}(X)_{K_X+N\geq 0}+\sum_{i=1}^nC_i+\sum_{j\in\Lambda}R_j,$$
    where each $R_j$ is a $(K_X+N)$-negative extremal ray that is relative ample at infinity with respect to $(X,N)$. Since $(K_X+L)\cdot C_i\geq 0$ and $P\cdot C_i=0$ for each $i$, $(K_X+N)\cdot C_i\geq 0$ for each $i$. Therefore,
        $$\overline{NE}(X)=\overline{NE}(X)_{K_X+N\geq 0}+\sum_{j\in\Lambda}R_j,$$
        and $R_j\not\subset\Supp N$ for any $j$. Therefore, $K_X\cdot R_j<0$ for each $j$. We let
    $$\Lambda_0:=\{j\in\Lambda\mid (K_X+L)\cdot R_j<0\}.$$
Since $P$ is nef,
    $$\overline{NE}(X)=\overline{NE}(X)_{K_X+L\geq 0}+\sum_{j\in\Lambda_0}R_j.$$  

        If $\Lambda_0=\emptyset$, then $K_X+L$ is nef, and the lemma follows from Lemma \ref{lem: nef case}. Otherwise, we pick $j_0\in\Lambda_0$ and let $\pi: X\rightarrow\bar X$ be the contraction of the $K_X$-negative extremal ray $R_{j_0}$. Since $\pi$ is a step of a $(K_X+L)$-MMP and $K_X+L$ is pseudo-effective, $\pi$ is birational, and $\bar X$ is smooth. Let $\bar L:=\pi_*L$, then $\bar L$ and $K_{\bar X}+\bar L$ are pseudo-effective, and $\rho(\bar X)=\rho(X)-1$. By Theorem \hyperlink{thm: A}{A}$_{\rho-1}$, $K_{\bar X}+\bar L$ is numerically effective. Since $\pi$ is $(K_X+L)$-negative,
        $$K_X+L=\pi^*(K_{\bar X}+\bar L)+E$$
        for some $E\geq 0$. Therefore, $K_X+L$ is numerically effective.
\end{proof}

\begin{proof}[Proof of Theorem \hyperlink{thm: A}{A}]
    We apply induction on $\rho$. The $\rho=1$ case is obvious as pseudo-effective and semi-ample are equivalent.

    Assume that $\rho\geq 2$. We let $P'$ and $N'$ be the positive and negative parts of the Zariski decomposition of $K_X+L$ respectively and let $Q:=N\wedge N'$. Then $Q\geq 0$, $P,\bar N:=N-Q$ are the positive and negative parts of the Zariski decomposition of $\bar L:=L-Q$ respectively, $P',\bar N':=N'-Q$ are the positive and negative parts of the Zariski decomposition of $K_X+\bar L$ respectively, and $\Supp\bar N\cap\Supp\bar N'=\emptyset$. Therefore, for any irreducible component $C$ of $\Supp\bar N$,
    $$(K_X+\bar L)\cdot C=(P'+\bar N')\cdot C\geq 0.$$
    By Lemma \ref{lem: nef on n case}, $K_X+\bar L$ is numerically effective. Thus $K_X+L=K_X+\bar L+Q$ is numerically effective.
\end{proof}

\section{Proof of the main theorems}

\begin{proof}[Proof of Theorem \ref{thm: main theorem}]
    Let $h: \tilde X\rightarrow X$ be the minimal resolution of $X$. Then $K_{\tilde X}+B=h^*K_X$ for some $B\geq 0$. Let $\tilde L:=B+h^*L$, then $\tilde L$ is pseudo-effective and $K_{\tilde X}+\tilde L=h^*(K_X+L)$ is pseudo-effective. By Theorem \hyperlink{thm: A}{A}, $K_{\tilde X}+\tilde L\equiv \tilde D\geq 0$ for some $\Qq$-divisor $D$ on $\tilde X$. Let $D:=h_*\tilde D$. Since $X$ is $\Qq$-factorial, by Lemma \ref{lem: num 0 pushforward 0}, $K_X+L\equiv D\geq 0$.
\end{proof}

\begin{proof}[Proof of Theorem \ref{thm: anti pe is ne}]
  Let $L:=-2K_X$. Then $L$ and $K_X+L$ are pseudo-effective. By Theorem \ref{thm: main theorem}, $K_X+L=-K_X$ is numerically effective.
\end{proof}

\begin{rem}
Theorem \ref{thm: main theorem} does not hold for $\Rr$-divisors. In fact, \cite[Example 3.19]{HL22} constructed a nef $\Rr$-divisor $M$ on $X:=E\times E$ and irreducible curves $C_n$ such that $0<M\cdot C_n<\frac{1}{n}$ for any positive integer $n$, where $E$ is a smooth elliptic curve. It is clear that $K_X+M=M$ is pseudo-effective but not numerically effective. In fact, if $K_X+M$ is numerically effective, then $K_X+M=M\equiv\sum_{i=1}^m a_iD_i$, where each $D_i$ is an irreducible curve and $a_i>0$ for each $i$. Thus for any $n\gg 0$, $C_n\not=D_i$ for any $i$. Since $M\cdot C_n>0$, $D_i\cdot C_n>0$ for some $i$, so $M\cdot C_n\geq\min\{a_i\mid 1\leq i\leq m\}$ for any $n$. This is not possible.

Of course, Theorem \ref{thm: main theorem} may still hold for $\Rr$-divisors after adding some restrictions, similar to the ``NQC" condition for nef $\Rr$-divisors (cf. \cite[Definition 2.15]{HL22}). We leave it to the readers.
\end{rem}


\begin{thebibliography}{99}

\bibitem[Amb03]{Amb03} F. Ambro, \textit{Quasi-log varieties}, Tr. Mat. Inst. Steklova \textbf{240} (2003), Biratsion. Geom. Linein. Sist. Konechno Porozhdennye Algebry, 220--239; translation in Proc. Steklov Inst. Math. \textbf{240} (2003), no. 1, 214--233.

\bibitem[BS23]{BS23} F. Bernasconi and L. Stigant, \textit{Semiampleness for Calabi-Yau surfaces in positive and mixed characteristic}, Nagoya Math. J., \textbf{250} (2023), 365--384.

\bibitem[BCHM10]{BCHM10}
C. Birkar, P. Cascini, C. D. Hacon and J. M\textsuperscript{c}Kernan, \textit{Existence of minimal models for varieties of log general type}, J. Amer. Math. Soc. \textbf{23} (2010), no. 2, 405--468.

\bibitem[BH14]{BH14} C. Birkar and Z. Hu, \textit{Polarized pairs, log minimal models, and Zariski decompositions}, Nagoya Math. J. \textbf{215} (2014), 203--224.

\bibitem[BZ16]{BZ16} C. Birkar and D.-Q. Zhang, \textit{Effectivity of Iitaka fibrations and pluricanonical systems of polarized pairs}, Publ. Math. IHÉS \textbf{123} (2016), 283--331.

\bibitem[Fon24]{Fon24} C. Fontanari, \textit{A question about generalized non-vanishing}, Rend. Circ. Mat. Palermo, II. Ser (2024).

\bibitem[Fuj17]{Fuj17} O. Fujino, \textit{Foundations of the minimal model program}, MSJ Memoirs \textbf{35} (2017), Mathematical Society of Japan, Tokyo.

\bibitem[HL22]{HL22} J.~Han and Z.~Li, \textit{Weak Zariski decompositions and log terminal models for generalized polarized pairs}, Math. Z. \textbf{302} (2022), 707--741.

\bibitem[HLS19]{HLS19} J. Han, J. Liu, and V. V. Shokurov, \textit{ACC for minimal log discrepancies of exceptional singularities}, arXiv:1903.04338.

\bibitem[HL20]{HL20} J. Han and W. Liu, \textit{On numerical nonvanishing for generalized log canonical pairs}, Doc. Math. \textbf{25} (2020), 93--123.

\bibitem[KMM94]{KMM94} S. Keel, K. Matsuki, and J. M\textsuperscript{c}Kernan, \textit{Log abundance theorem for threefolds}, Duke Math. J. \textbf{75} (1994), 99--119.

\bibitem[KM98]{KM98} J. Koll\'{a}r and S. Mori, \textit{Birational geometry of algebraic varieties}, Cambridge Tracts in Math. \textbf{134} (1998), Cambridge Univ. Press.

\bibitem[LMPTX23]{LMPTX23} V. Lazi\'c, S. Matsumura, T. Peternell, N. Tsakanikas, and Z. Xie, \textit{The Nonvanishing problem for varieties with nef anticanonical bundle}, Doc. Math. \textbf{28} (2023), no. 6, 1393--1440.

\bibitem[LP20a]{LP20a} V. Lazi\'c and T. Peternell, \textit{On Generalised Abundance, I}, Publ. Res. Inst. Math. Sci. \textbf{56} (2020), no. 2, 353--389.

\bibitem[LP20b]{LP20b} V. Lazi\'c and T. Peternell, \textit{On generalised abundance, II}, Peking Math. J. \textbf{3} (2020), 1--46.



\bibitem[Mul23]{Mul23} N. M\"uller, \textit{Two Non-Vanishing results concerning the Anti-Canonical Bundle}, arXiv:2305.19109.





\end{thebibliography}
\end{document}